\documentclass[12pt]{article}
\usepackage{amssymb}
\usepackage{graphicx}
\usepackage{subcaption}
\usepackage{amsfonts}
\usepackage{latexsym}
\usepackage{amsthm}
\usepackage{amsmath}

\usepackage{amssymb, amscd}

\usepackage{url}
\usepackage[T1]{fontenc}
\usepackage{color}
\usepackage{tabu}

\usepackage{bbm}

\usepackage{amsmath,amssymb,amscd,}
\usepackage[all]{xy}
\usepackage{color}

\numberwithin{equation}{section}

\newtheorem{conj}{Conjecture}[section]
\newtheorem{thm}[conj]{Theorem}

\newtheorem{defin}[conj]{Definition}
\newtheorem{prop}[conj]{Proposition}

\newtheorem{lema}[conj]{Lemma}

\begin{document}

 \title{An elementary approach to  local combinatorial formulae for the Euler class of a PL spherical fiber bundle}

\author{ Gaiane Panina\\
{\small PDMI} \\[-2mm]
}

\maketitle
\begin{abstract}\noindent
We present an elementary approach to   local combinatorial formulae for the Euler class of a fiber-oriented triangulated spherical fiber bundle.
 The approach is based on \textit{sections averaging} technique and very basic knowledge of simplicial (co)homology theory.  Our formulae are  close
relatives of  those by N. Mn\"{e}v.
\end{abstract}

\section{Introduction}
\label{sec:intro}
Recently N. Mn\"{e}v presented a local combinatorial formula for the Euler class of a fiber-oriented $ n$-dimensional PL spherical  fiber bundle in terms of
combinatorial  Hodge theory twisting cochains in
Guy Hirsch's homology model \cite{M2}.
The present paper introduces an elementary approach to local combinatorial formulae.

Let $S^n\rightarrow E \xrightarrow[\text{}]{\pi} B$ be a locally trivial fiber bundle whose fiber is the oriented sphere $S^n$.
We  assume that $B$ and $E$ are triangulated in such a way that $\pi$ is a \textit{simplicial map}, that is, $\pi$ maps simplices to simplices, linearly on each of the simplices.
Besides, we assume that  $\pi^{-1}(x)$   $ \forall x \in B$ is a \textit{ combinatorial manifold}. This means in particular that the dual cell decomposition is well-defined. 

A local combinatorial formula for the rational Euler class $e(E \xrightarrow[\text{}]{\pi} B)$ is 
 an algorithm which associates a rational number to every oriented $(n+1)$-dimensional simplex $\sigma^{n+1} $ from the base.  The output of the algorithm is a cochain
$\mathcal{E}$  representing the Euler class. The value of the cochain at a simplex $\sigma^{n+1}$ depends only on the restriction of the bundle $\pi^{-1}(\sigma^{n+1}) \rightarrow \sigma^{n+1}$.

We present two formulae. The first one coincides with that of N. Mn\"{e}v, Section \ref{SecMain1}. The second  one presents some shortcul on the last step of the algorithm, see Section \ref{SecMain}.

One can alter these algorithms even further (we explain how in Section \ref{sectionconcl}) and thus obtain various different local combinatorial formulae.

In Section \ref{sectionconcl} we explain that for  $n=1$ our second formula (Section \ref{SecMain}) coincides with the Igusa-Mn\"{e}v-Sharygin combinatorial formula for circle bundle, see \cite{I} and \cite{M1}.

\bigskip

Here is the\textbf{ leading  idea  in short} (it traces back to M. Kazarian's\textit{ multisections} \cite{Kaz}, used also in \cite{G}).
 There exists the following way to compute the Euler class:
 
 \begin{prop}\label{prop-Euler} \cite{Miln},\cite{Fom}

Let $E \rightarrow B$ be a fiber-oriented spherical bundle over a triangulated base $B$. Assume that  its partial section $s$   is defined on the $n$-skeleton of $B$  .
Fix an orientation on each of the $(n+1)$-dimensional simplices $\sigma^{n+1}\in B$, and  set $\mathcal{E}(s,\sigma^{n+1})$ to be the degree of the map $s:\partial \sigma^{n+1} \rightarrow S^n$.
We claim that the integer cochain $\mathcal{E}(s, \sigma^{n+1})$ is a cocycle and  represents the Euler class $e(E \xrightarrow[\text{}]{\pi} B)$.
\end{prop}
For bundles arising from  vector bundle with the base $ B$ of rank $(n+1)$, this proposition is proven in J. Milnor and Stasheff  \cite{Miln}. The general case is treated in A. Fomenko and D. Fuchs  \cite{Fom}, 23.5. \qed

\medskip

Next, if we have several partial sections $s_1,...,s_r$, set $\mathcal{E}_i(\sigma^{n+1})=\mathcal{E}(s_i,\sigma^{n+1})$. The average  $\frac{1}{r}\sum_i \mathcal{E}_i(\sigma^{n+1})$ is clearly a rational cochain representing the Euler class, but not necessarily integer. Furthermore,  if we have partial sections $s_1,...,s_M$ and $s_1',...,s_N'$,
 the following average  of the cochains  $$\frac{1}{M-N}\Big(\sum_i \mathcal{E}_i-\sum_j \mathcal{E}_j'\Big)$$
also represents the rational Euler class provided that $M \neq N$.

The combinatorics of the triangulation of the bundle suggests a way of fixing a collection of partial sections.
We create these sections stepwise, starting by $0$-skeleta, and extending them further. Our toolbox contains \textit{ harmonic chains } and \textit{winding numbers}. 
\medskip

\textbf{Reader's guide:} For a shortcut to the first main result (Section \ref{SecMain1}) which coincides with Mn\"{e}v's formula, it suffices to read  Section \ref{subsec21},  and  Definition \ref{DefHarm}.
 For a shortcut to the second main result (Section \ref{SecMain}) one needs also Definitions \ref{DefW},
and \ref{DefW2}.

\textbf{Acknowledgement:} The author is grateful to Alexander Gaifullin, Nikolai Mn\"{e}v and Ivan Panin  for useful comments.

\section{Toolbox}

Denote the associated simplicial complexes of the base and of the total space by $K_B$ and $K_E$ respectively.

\subsection{Combinatorics of preimage of a simplex}\label{subsec21}

Fix an oriented simplex $\sigma^{n+1} \in K_B$ in the base. Denote its vertices by $v_0,...,v_{n+1}$ in such a way that the order agrees with the orientation of $\sigma^{n+1}$.
It is  instructive to think that vertices are colored, say, $v_0$ is red, $v_1$ is blue, $v_2$ is green, etc.

Let us analyse the $2n+1$-dimensional simplices lying in the preimage $\pi^{-1}(\sigma^{n+1})$. Each such simplex $\Delta^{2n+1}\in \pi^{-1}(\sigma^{n+1})$ has $2n+2$ vertices colored   according to the colors of their projections. Clearly, for each color $i$ there is a vertex of $\Delta^{2n+1}$ colored by $i$.

We
 say that  $\Delta^{2n+1} \in \mathcal{A}_i(\sigma^{n+1})$ if $$|Vert(\Delta^{2n+1})\cap \pi^{-1}(v_i)|=n+1,$$ that is, $\Delta^{2n+1}$ has $n+1$ vertices of   color $i$, whereas other colors appear only once.  Note that the classes $\mathcal{A}_i$ do not cover all the simplices in the preimage.


Now analyse $\pi^{-1}(v_i)$. It is a triangulated oriented sphere $S^n$. The maximal simplices of the triangulation  correspond to simplices from  $\mathcal{A}_i(\sigma^{n+1})$. Denote the dual  of the triangulation by $\Gamma_i$. It is a regular cell complex, and we  imagine it  to be colored in color $i$.

\begin{figure}[h]
\centering \includegraphics[width=10 cm]{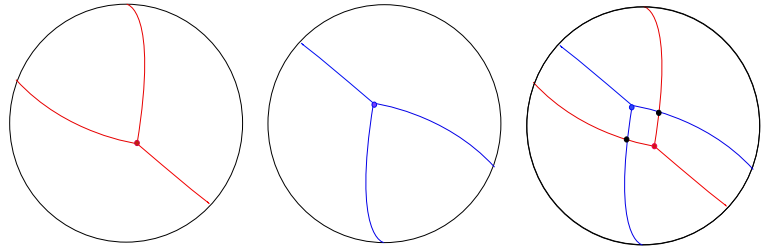}
\caption{Complexes $\Gamma_1, \ \Gamma_2$, and their superposition $\Gamma_{12}$.}\label{fig1}
\end{figure}

Now  take a point $x$ lying strictly inside an  edge, say, in $v_0v_1$.  Its preimage  $\pi^{-1}(x)$ (independently on the choice of $x$) is a tiling of $S^n$, whose dual complex is the \textit{superposition} of $\Gamma_0$ and $\Gamma_1$, that is, the mutual tiling  generated by  $\Gamma_0$ and $\Gamma_1$, see Fig. \ref{fig1}.

Generalizing, we conclude that:
\begin{prop}\begin{enumerate}
              \item For any point  $x\in \sigma^{n+1}$, the combinatorics of the preimage $\pi^{-1}(x)$ depends only on the face  $F$ of $\sigma^{n+1}  $ that contains $x$ as an inner point.

              \item Let $\{v_{i_1},...,v_{i_k}\}$  be the set of vertices of $F$. Then the
dual complex $\Gamma_{i_1,...,i_k}$ of the preimage $\pi^{-1}(x)$ is the superposition of  $\Gamma_i$ where $i$ ranges over the set $\{i_1,...,i_k\}$.
              \item
In particular, for an inner point of $\sigma^{n+1}$, we have the superposition of $(n+2)$ colored complexes  $\Gamma_0,...,\Gamma_{n+2}$.
Elimination of one (or several) of these complexes corresponds to moving $x$ to some face.\qed
            \end{enumerate}
\end{prop}

Some of the  vertices of $\Gamma_{i_1,...,i_k}$  are the vertices of $\Gamma_i$ for some $i$. Let us say that these vertices \textit{are colod by color} $i$.
The other (uncolored vertices) arise as intersections of the cells of $\Gamma_i$ for different $i$, see Fig. \ref{fig1}  for an illustration.

\begin{lema}\label{LemmaNCells} Let $V$ be a   vertex of $\Gamma_{i}$. Then
\begin{enumerate}
  \item
$V$ is also a vertex of the complex $\Gamma_I$  provided that $i\in I$, and
  \item
 $V$ has exactly $n+1$ incident top-dimensional cells of $\Gamma_I$. \qed.
\end{enumerate}
\end{lema}

\bigskip

\subsection{S-chains, patches, and s-patches}
Assume there is a  regular cell complex $\Gamma$ which is a combinatorial sphere $S^n$. Denote by $C_{k+1}(\Gamma, \mathbb{Q})$ the rational chain group of $\Gamma$.  Let $k<n$, and let $a\in C_k(\Gamma, \mathbb{Q})$ be  a closed $k$-dimensional chain.
A chain $c\in C_{k+1}(\Gamma, \mathbb{Q})$ is called \textit{a patch}   of $a$ if $\partial c=a$.

Assume  there is a map  from the oriented  $k$-sphere to the $k$-skeleton $$s:S^k \rightarrow Skel_k(\Gamma).$$

It yields a map $$S^k \rightarrow Skel_k(\Gamma)/Skel_{k-1}(\Gamma).$$

The factor space $Skel_k(\Gamma)/Skel_{k-1}(\Gamma)$ is a wedge of $k$-spheres, where the spheres correspond to $k$-cells of $\Gamma$.


We conclude  that $s$ gives  rise to some closed chain $|s|\in C_k(\Gamma, \mathbb{Z})$  which represents the pushforward of the fundamental class of $S^k$.

Assume now there is a map  $$\tilde{s}:D^{k+1}\rightarrow Skel_{k+1}(\Gamma)$$  such that the boundary of $D^{k+1}$ is mapped to $Skel_k(\Gamma)$. Similarly, the map $\tilde{s}$ gives  rise to a $k+1$-chain $|\tilde{s}|\in C_{k+1}(\Gamma, \mathbb{Z}).$  Indeed,  in this case $s$
yields a map $$D^{k+1}/\partial D^{k+1}\rightarrow Skel_{k+1}(\Gamma)/Skel_{k}(\Gamma).$$

\bigskip

If a chain $a$ equals $|s|$  for some $s$, we say that $s$ \textit{supports} $a$. In this case we call $a$ an \textit{s-chain}.


Assume that  $a=|s|$ for  some $s:S^k \rightarrow Skel_k(\Gamma)$.

Its patch  $c$ is called an \textit{s-patch  of $a$} if $c=|\tilde{s}|$ for an extension $\tilde{s}:D^{k+1}\rightarrow Skel_{k+1}(\Gamma)$ of the map $s$.

\bigskip
\textbf{
Example.}  For $k=0$ let $a$ be an ordered pair of vertices $V_0,V_1$ interpreted as an s-chain.  Each path on $Skel_1$ connecting the vertices is an s-patch of $a$. However not every patch of $a$ is an s-patch. For instance the chain that averages two different paths is not an s-patch.

\bigskip

\begin{lema}\label{LemmaAver}

1.  If $c=|\tilde{s}|$ is a chain supported by $\tilde{s}:D^{k+1} \rightarrow Skel_{k+1}(\Gamma)$, and \newline $a$ is a chain supported by $$s=\tilde{s}_{|_{\partial  D^k}}:\partial  D^k \rightarrow Skel_k(\Gamma)$$
then $\partial c =a.$

 2. Each rational closed $k$-chain is a rational linear combination of s-chains  supported by some $s_i:S^{k} \rightarrow Skel_{k}(\Gamma)$.

3.   Let $k<n$. Let $a=|s|$ be an s-chain supported by some $s:S^{k} \rightarrow Skel_{k}$. Let $b$ be a rational chain   such that $\partial b =a.$  Then
$b$  is\textbf{ an average }of chains supported by extensions of the map $s$.

More precisely, there exist

 $$\tilde{s}_i: D^{k+1} \rightarrow Skel_{k+1}(\Gamma), \ \ i=1,...,M$$  and
 $$\tilde{s}'_j: D^{k+1} \rightarrow Skel_{k+1}(\Gamma), j=1,...,N$$

 such that  the restrictions of all the  $\tilde{s}_i$ and  $\tilde{s}'_j$  to $\partial D^{k+1}$ equal $s$,
 and
 $$b=\frac{1}{M-N}\Big(\sum_{i=1}^M |\tilde{s}_i|-\sum_{j=1}^N |\tilde{s}'_j|\Big).$$
\end{lema}
Proof.

(1) is trivial.

(2) Let $k<n$. Let $\tau$ range over all the $k+1$-cells of the complex $\Gamma$.  Then $\{\partial \tau \}$  generate the space of closed $k$-chains.  Indeed,
for any closed chain $p \in C_k(\Gamma)$, for some $r\in C_{k+1}(\Gamma)$ we have $p=\partial r=\partial (\sum  a_i \tau_i)= \sum a_i \partial \tau_i$.

\bigskip

If $k=n$ the proof is even simpler:  all closed chains differ by a multiple.

\bigskip

(3) Let $k<n-1$. Take any extension of $s$ to a map $\tilde{s}$ from the disk $D^{k+1}$ to the $k$-skeleton, that is, $\tilde{s}_{|_{\partial D}}=s$.  Then $b-|\tilde{s}|$ is a closed rational $k+1$-chain and hence a rational linear combination of $ \partial \tau_i$ where $\tau_i$ ranges over $(k+2)$-cells. Let us prove that for each $(k+2)$-cell $\tau$, the chain $\partial (\tau)$ is a difference of some $|\tilde{s}|$ and  $|\tilde{s}'|$   such that  the restrictions of   $\tilde{s}$ and  $\tilde{s}'$  to the boundary of the disk equal $s$.

Fix any subcomplex of $\partial \tau$  called $Equator$ which is homeomorphic to $S^{k-1}$.  $Equator$  splits $\partial \tau$ into two  hemispheres $B_1$ and $B_2$.  More precisely,  denote by $B_1$ and $B_2$ the chains such that $B_1+B_2=\partial \tau$. Take a  map of the cylinder $S^k\times[0,1]$ to the $k$-skeleton such that the restriction of the map on $S^k \times 0$ equals $s$, and the restriction on $S^k \times 1$  is a homeomorphism to the $Equator$. Such a map exists since any two maps from the $k$-sphere to the $k+1$-skeleton are homotopic. We treat the cylinder as a disk $D^{k+1}$ with an eliminated smaller $k+1$-disk. Extend the map to the entire $D^k$ in two ways, using $B_1$ and  $B_2$. Denote the extensions by  $\tilde{s}$ and $\tilde{s}'$.
Now we have
 $\partial \tau= |\tilde{s}|-|\tilde{s}'|$.

 \bigskip

 For  $k=n-1$ the proof is again simpler. In this case two rational patches of $a$ differ by a rational multiple of the chain representing the fundamental class of $S^n$.
 Each integer patch of $a$ is an $s$-patch.
 \qed

\medskip

In this lemma we  call the entries $s_i$ (resp., entries $s'_j$) \textit{ positive} (resp. \textit{negative})  entries of the averaging.

\newpage
\subsection{Harmonic extension}
Harmonic extension gives a way  to canonically reverse the boundary operator.

We borrow the following basic knowledge about harmonic chains from  \cite{E}.
Given a  regular cell complex $K$, we assume that $C_{k}(K, \mathbb{Q})$ is endowed by a scalar product such that the collection of all the $k$-cells correspond to an orthonormal basis.  Let $\partial$ denote the standard boundary operator, and let $\partial^*$ be its adjoint operator  related to the scalar product. Define the\textit{ discrete Laplacian} as $\Delta= \partial^*\partial + \partial \partial^*$.
 The kernel of $\Delta$ consists of some closed chains (called \textit{harmonic chains}). It is known that each homology class is uniquely representable by a harmonic chain.

If $\Delta:C_k(K, \mathbb{Q})\rightarrow C_k(K, , \mathbb{Q})$ is invertible (i.e., if $H_k(K,  \mathbb{Q})$ vanishes), there exists the inverse operator $\Delta^{-1}:C_k(K, \mathbb{Q})\rightarrow C_k(K, \mathbb{Q})$ called \textit{Green operator}.

\bigskip

Assume we have a  regular cell complex $\Gamma$ which is a combinatorial sphere $S^n$.
Take  the\textit{ cone over} $Skel_k(\Gamma)$. It is a  cell complex defined as   $Cone(Skel_k(\Gamma))=Skel_k(\Gamma)* O$, where  $*$ denotes the join operation, and $O$ is the apex of the cone.  Patch the cone to the combinatorial sphere $\Gamma$  via the natural inclusion map $Skel_k \subset Cone(Skel_k)$. Also patch a disk $D^{n+1}$ to $\Gamma$ along a homemorphism between $\Gamma$ and $\partial D^{n+1}$.  We obtain some cell complex  and denote it by $PC_k(\Gamma)$.

Let $k<n-1$, and let $a\in C_k(\Gamma, \mathbb{Q})$ be  a closed  chain. Consider a (uniquely defined) $(k+1)$-chain $\overline{a}\in C_{k+1}(Cone(Skel_k(\Gamma))$  with the property $\partial \overline{a}=-a$. In simple words, $\overline{a}$ assigns the coefficient $-a(\sigma^k)$ to the cell $\sigma^k*O$.  Now we extend the chain by keeping its coefficients on the cells of type $\sigma^k*O$ and assigning some coefficients to $(k+1)$-cells of $\Gamma$.  The chain $\overline{a}$ extends in many ways to a closed chain on $PC_k(\Gamma)$. All closed extensions represent one and the same homology class of $PC_k(\Gamma)$. Denote the unique harmonic (with respect to $PC_k(\Gamma)$) representative   of the homology class  by $h(a)$.

 Decompose  $h(a)$  as  $h(a)=\overline{a}+\mathcal{H}(a)$. We have necessarily $\mathcal{H}(a)\in C_{k+1}(\Gamma)$. Since $\Delta h(a)=0$, we have
$\Delta \overline{a} +\Delta \mathcal{H}(a)=0$, and therefore $\Delta \mathcal{H}(a) =\partial^* a$.

Since $h(a)$ is a closed chain, $\partial h(a)= 0$,  so $\partial \overline{a}+\partial \mathcal{H}(a)=0$, hence $\partial \mathcal{H}(a) =a$.

Thus we arrive at the following:

\begin{defin} \label{DefHarm} Let $k\leq n-1$, and let $a\in C_k(\Gamma, \mathbb{Q})$ be  a closed $k$-dimensional chain. The \textit{harmonic extension} of $a$  is the
 chain $\mathcal{H}(a)\in C_{k+1}(\Gamma, \mathbb{Q})$  defined as $ \mathcal{H}(a)= \Delta^{-1}(\partial^* a)$.
 \end{defin}

\begin{lema}\label{LemmaHarmSpatch}
\begin{enumerate}

  \item For the harmonic extension we have
$$\partial \mathcal{H}(a)= a.$$
\item The harmonic extension  of an s-chain $a$ is representable as

$$\mathcal{H}(a)=\frac{1}{N-M}\Big(\sum_{i=1}^N  c_i-\sum_{j=1}^M  c'_j\Big),$$$ \text{  where }  M \neq \text{$N$  are some natural numbers, }   \ \partial c_i=\partial c_j' =a,$
and each of the chains $c_i$ and $c_j'$  is an s-patch  of $a$.\qed
\end{enumerate}

\end{lema}

\medskip

\textbf{Important remarks:}  1. Harmonic extension depends on the complex $\Gamma$, not only on the chain $a$. For instance, if one subdivides $\Gamma$ but keeps $a$, the harmonic extension changes.

2. Although harmonic extension makes sense for all $k\leq n-1$, for the sake of a shortcut, for $n-1$-dimensional chains we use  also a  different technique of reversal of the boundary operator, see below.

\medskip

\subsection*{Winding number. Dimension-$n$ extension}  Throughout the section $\Gamma$ is a regular cell complex which is a combinatorial sphere  $S^n$. Let $\Sigma \in C_{n-1}(\Gamma, \mathbb{Q})$  be a closed chain, let $x,y\in S^n\setminus Skel_{n-1}(\Gamma)$  be some points.

\begin{defin}\label{DefW}
The winding number $\mathcal{W}(x,y,\Sigma)$ is defined as the algebraic number of intersections of $\Sigma$ and $[x,y]$, where $[x,y]$ is a smooth (or a piecewise linear) oriented path  from $x$ to $y$ which intersects the $(n-1)$-skeleton of $\Gamma$ transversally. If the endpoints $x$ and $y$ are fixed, the winding number is independent on the choice of the path $[xy]$.  Borrowing notation from intersection theory, and treating $[x,y]$ as a $1$-chain, we write
$$\mathcal{W}(x,y,\Sigma)=\Sigma \smile [x,y].$$
\end{defin}


\bigskip

\begin{figure}[h]
\centering \includegraphics[width=5 cm]{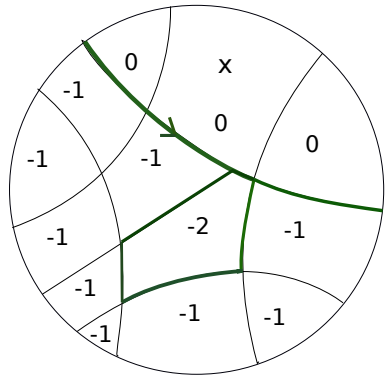}
\caption{Winding number for the chain depicted as the green curve.}\label{fig2}
\end{figure}

\textbf{Examples.}\begin{enumerate}
                    \item If $\Sigma$ arises as a pushforward of the fundamental class of some $\psi:S^1\rightarrow S^{2}$, and one imagines $S^2$ as $\mathbb{R}^2$  compactified by the point $\infty=y$, then $\mathcal{W}(x,y,\Sigma)$ is exactly the winding number of $\Sigma$ around the point $x$, see Fig. \ref{fig2} for an illustration.
                    \item If $\Sigma$ is an average of cycles represented by submanifolds, then $\mathcal{W}(x,y,\Sigma)$ is the average of winding numbers.
                  \end{enumerate}

\bigskip

\begin{lema}

 Let $\Sigma\in C_{n-1}(\Gamma, \mathbb{Q})$  be a closed chain.
\begin{enumerate}
  \item For each cell $\sigma_1^n\in \Gamma$ there exists a unique chain $\mathcal{C}=\mathcal{C}_{\Sigma,\sigma^n_1}\in C_{n}(\Gamma, \mathbb{Q})$ whose coefficient $\mathcal{C}(\sigma_1^n)$ at $\sigma_1^n$ vanishes, and
$\partial \mathcal{C}=\Sigma$.
  \item The values of $\mathcal{C}$ on other cells are the winding numbers:
For all $\sigma_2^n \in \Gamma$,
$$\mathcal{C}(\sigma_2^n)= \mathcal{W}(x,y,\Sigma),$$
where $x\in \sigma_1, \ y \in \sigma_2$.
\item If $\Sigma$ is an s-chain, then $\mathcal{C}$ is an s-patch.
   \qed

\end{enumerate}
\end{lema}

\bigskip

\begin{defin}\label{DefW2}
 Let $\Sigma \in C_{n-1}(\Gamma, \mathbb{Q})$  be a closed chain, let  $y\in S^n\setminus Skel_{n-1}(\Gamma)$ be some point, and $V$ is a vertex of \ $\Gamma$. Assume also that $x$ is incident to exactly $n+1$ top-dimensional cells of $\Gamma$.

The winding number is defined as the average $$\mathcal{W}(V,y,\Sigma)=\frac{1}{n+1}\sum\mathcal{W}(x_i,y,\Sigma),$$ where the sum ranges over points $x_i$, one point from each $n$-cell incident to $V$. From now on for fixed $V$ and $\Sigma$ we treat $\mathcal{W}(V, \Sigma)$ as an $n$-chain, the linear combination of $n$-cells, such that the coefficient of a cell equals the winding number of an inner point of the cell.

\end{defin}

The following is obvious:
\begin{lema}\label{LemmaDecompN-ext}
\begin{enumerate}
\item
$\partial \mathcal{W}(V,y,\Sigma)=\Sigma$.
\item The operator $ \mathcal{W}(V,y,\Sigma)$  is  linear with respect $\Sigma$.
 \item  If $\Sigma$ is an s-chain, then $ \mathcal{W}(V,\Sigma)$  is an average of its s-patches.
\qed

\end{enumerate}
\end{lema}

\bigskip

\newpage

\section{The first formula}\label{SecMain1}

\begin{enumerate}
  \item Fix a simplex $\sigma^{n+1}$ in the base $B$. Denote its  vertices by $v_0,...,v_{n+1}$ and consider the complex $\Gamma_{01...n+1}=\Gamma_{01...n+1}(\sigma^{n+1})$. Take all existing $(n+2)$-tuples  $V_0,...,V_{n+1}$ of its   vertices, one of each color, that is, $V_i\in \pi^{-1}(v_i)$  is a vertex of $\Gamma_i$. For each of the tuples do the following.
  \item   For each oriented edge $ij$, take the $0$-chain $-V_i+V_j\in C_0(\Gamma_{ij}, \mathbb{Q})$ and compute its {harmonic extension}  related to the complex $\Gamma_{ij}$:  $$\Sigma_{ij}=\mathcal{H}(-V_i+V_j).$$

  \item For each oriented triple of vertices $ijk$  compute the harmonic extension related to the complex $\Gamma_{ijk}$: $$\Sigma_{ijk}=\mathcal{H}( \Sigma_{jk}-\Sigma_{ik}+\Sigma_{ij}).$$

      \bigskip

\item   Proceed in the same way  for $k\leq n$. Consider a $(k+1)$-tuple of the vertices  of  the simplex $\sigma^{n+1}$, say, $v_0,...,v_{k}$.
Assuming that the chains  $\Sigma_{0...\hat{i}...k}$ are defined on the previous step,
    compute the harmonic extension  related to the complex $\Gamma_{0...k}$: $$\Sigma_{0...k}=\mathcal{H}\Big(\sum_{i=0}^{k}(-1)^i\Sigma_{0...\hat{i}...k}\Big)\in C_k(\Gamma_{0...k}).$$

  \item  Eventually we arrive at  a collection of chains $\Sigma_{012...\hat{i}...n+1}$. The chain
  $$\sum_{i=0}^{n+1}(-1)^i\Sigma_{012...\hat{i}...n+1}\in C_{n}(\Gamma_{01...n+1})$$
  is a closed $n$-chain, and thus represents some rational number
$\mathfrak{e}=\mathfrak{e}(\sigma^{n+1}).$

\end{enumerate}

\begin{thm}\label{ThmMain1}
\begin{enumerate}
 
\item The number $\mathfrak{e}(\sigma^{n+1})$
  depends on the combinatorics of $\pi^{-1}(\sigma^{n+1})$ and the choice of $V_1,...,V_n$ only.
   \item A change of orientation of $\sigma^{n+1}$ changes the sign of $\mathfrak{e}(V_0,...,V_{n+1};\Gamma_{01...n+1})$.
  \item The function that assigns to each $(n+1)$-dimensional simplex the $\sigma^{n+1}$ the average of $\mathfrak{e}(V_0,...,V_{n+1};\Gamma_{01...n+1})$ over all $(n+2)$-tuples of vertices   $(V_0,...,V_{n+1})$, one vertex from each color,
      $$\mathcal{E}(\sigma^{n+1})=\frac{1}{|Vert(\Gamma_0)|...|Vert(\Gamma_{n+1})|}\ \ \sum_{\forall \ i \ V_i\  \in \  Vert(\Gamma_i)} \mathfrak{c}(V_0,...,V_{n+1};\Gamma_{01...n+1})$$
       is a closed cochain that represents the rational Euler class of the fiber bundle.
       \item This cochain coincides with the cochain obtained in Mn\"{e}v  \cite{M2}, Theorem 8.
\end{enumerate}
\end{thm}

\section{The second formula: using a shortcut on the last step}\label{SecMain}

\begin{enumerate}
  \item Fix a simplex $\sigma^{n+1}$ in the base $B$. Denote its  vertices by $v_0,...,v_{n+1}$ and consider the complex $\Gamma_{01...n+1}=\Gamma_{01...n+1}(\sigma^{n+1})$. Take all existing $(n+2)$-tuples  $V_0,...,V_{n+1}$ of its   vertices, one of each color, that is, $V_i\in \pi^{-1}(v_i)$  is a vertex of $\Gamma_i$. For each of the tuples do the following.
  \item   For each oriented edge $ij$, take the $0$-chain $-V_i+V_j\in C_0(\Gamma_{ij}, \mathbb{Q})$ and compute its {harmonic extension}  related to the complex $\Gamma_{ij}$:  $$\Sigma_{ij}=\mathcal{H}(-V_i+V_j).$$

  \item For each oriented triple of vertices $ijk$  compute the harmonic extension related to the complex $\Gamma_{ijk}$: $$\Sigma_{ijk}=\mathcal{H}( \Sigma_{jk}-\Sigma_{ik}+\Sigma_{ij}).$$

      \bigskip

\item   Proceed in the same way  for $k<n$. Consider a $(k+1)$-tuple of the vertices  of  the simplex $\sigma^{n+1}$, say, $v_0,...,v_{k}$.
Assuming that the chains  $\Sigma_{0...\hat{i}...k}$ are defined on the previous step,
    compute the harmonic extension  related to the complex $\Gamma_{0...k}$: $$\Sigma_{0...k}=\mathcal{H}\Big(\sum_{i=0}^{k}(-1)^i\Sigma_{0...\hat{i}...k}\Big)\in C_k(\Gamma_{0...k}).$$

  \item  The last step is different: assuming that $\Sigma_{012...\hat{i}...\hat{j}...n+1}$ are already known, denote

$$\Sigma_{\hat{i}}=\Sigma_{012...\hat{i}...n+1}=\sum_{j\neq i} (-1)^{j+\frac{sign(j-i)+1}{2}}\Sigma_{012...\hat{i}...\hat{j}...n+1},$$

and set

$$\mathfrak{E}_i=\frac{1}{n+1}\sum_{k\neq i}\mathcal{W}(V_i,V_k,\Sigma_{\hat{i}})\ \in C_n(\Gamma_{0...n+1}, \mathbb{Q}).$$

\end{enumerate}

\begin{thm}\label{ThmMain}
\begin{enumerate}
  \item  The chain $\mathfrak{E}_i$    does not depend on $i$.
\item The chain $\mathfrak{E}_i$
is closed, and thus equals the fundamental class multiplied by some rational number $\mathfrak{e}=\mathfrak{e}(V_0,...,V_{n+1};\Gamma_{01...n+1})$  which  depends on the combinatorics of $\pi^{-1}(\sigma^{n+1})$ and the choice of $V_1,...,V_n$ only.
   \item A change of orientation of $\sigma^{n+1}$ changes the sign of $\mathfrak{e}(V_0,...,V_{n+1};\Gamma_{01...n+1})$.
  \item The function that assigns to each $(n+1)$-dimensional simplex the $\sigma^{n+1}$ the average of $\mathfrak{e}(V_0,...,V_{n+1};\Gamma_{01...n+1})$ over all $(n+2)$-tuples of vertices   $(V_0,...,V_{n+1})$, one vertex from each color,
      $$\mathcal{E}(\sigma^{n+1})=\frac{1}{|Vert(\Gamma_0)|...|Vert(\Gamma_{n+1})|}\ \ \sum_{\forall \ i \ V_i\  \in \  Vert(\Gamma_i)} \mathfrak{c}(V_0,...,V_{n+1};\Gamma_{01...n+1})$$

       is a closed cochain that represents the rational Euler class of the fiber bundle.
\end{enumerate}
\end{thm}

\section{Proof  of Theorems  \ref{ThmMain1} and  \ref{ThmMain}. Construction  of partial sections}

We shall construct a collection  of partial sections step by step, guided by combinatorics of the triangulated fiber bundle.
We start by partial sections over the zero skeleton of the complex $K_B$ and extend them stepwise to the skeleta of bigger dimensions  up to $Skel_n(K_B)$.

\medskip
We assume that an orientation is fixed for all the simplices of the complex $K_B$.
Let us start by fixing an oriented simplex $\sigma^{n+1} \in K_B$, enumerate its vertices consistently with the orientation,  and fix a local trivialization of the bundle over $\sigma^{n+1}$.

\medskip

\textbf{Step 0.}  Take one of the vertices of $\sigma^{n+1}$, say, $v_0$. Its preimage  $\pi^{-1}(v_0)$ is a triangulated sphere. Set a partial section $s$ over $v_0$ to be equal a vertex $V_0$ of the dual cell complex $\Gamma_0$ of the triangulation.

Repeat this for all the  vertices of the complex $K_B$. So, we have a collection of partial sections over $Skel_0(K_B)$.

In other words, for the simplex $\sigma^{n+1}$ a partial section fixes a collection $V_0,...,V_{n+1}$ of $n+2$  colored vertices in the complex $\Gamma_{012...n+1}$, one vertex for each color.

\medskip

\textbf{Step 1.} Take an edge of $\sigma^{n+1}$, say, $(v_iv_j)$, and consider the  complex $\Gamma_{ij}$. Remind that it equals the superposition of  $\Gamma_{i}$ and $\Gamma_{j}$. We have already defined a partial section over the vertices $v_i$ and $v_j$ by fixing $V_i$ and $V_j$.  These are vertices of the complex $\Gamma_{ij}$.
The chain
$$\Sigma_{ij} =\mathcal{H}(-V_i +V_j)\in C_1(\Gamma_{ij},\mathbb{Q})$$   is a (uniquely defined) $1$-chain    whose boundary is $-V_i +V_j$.
By Lemma \ref{LemmaHarmSpatch}, $\Sigma_{ij}$ is an average of a collection of  chains supported by  paths  from $V_i$ to $V_j$, that is, as an average of chains supported by partial sections over the edge $v_iv_j$.

Treating in similar way the rest of the edges and the rest of the simplices,
we arrive at a collection of partial sections over  $Skel_1(K_B)$.

\medskip

\textbf{
Step 2.}

Now our goal is to extend each of the partial sections  to  $Skel_2(K_B)$.  For each triple of the vertices $v_i,v_j,v_k$ of $\sigma^{n+1}$ consider the closed $1$-chain  $ \Sigma_{jk}-\Sigma_{ik}+\Sigma_{ij}$.
Here we tacitly assume that the simplex $v_iv_jv_k$ is oriented, and the order $ijk$ is consistent with the orientation.
Remind that $$\Sigma_{ijk}=\mathcal{H}( \Sigma_{jk}-\Sigma_{ik}+\Sigma_{ij}) \in C_2(\Gamma_{ijk}, \mathbb{Q}).$$

is the harmonic extension  taken with respect to $\Gamma_{ijk}$.

\begin{lema}\label{lemma41} The chain $\Sigma_{ijk}$ is an average of chains supported by partial sections of the bundle defined on the face $(v_iv_jv_k)$ of the simplex $\sigma^{n+1}$.
\end{lema}
Proof. We know already that $\Sigma_{ij}$, $\Sigma_{jk}$, and $\Sigma_{ki}$  are averages chains supported  by some paths. We may
assume that in the averaging expressions from Lemma \ref{LemmaAver}
 $$\frac{1}{M-N}\Big(\sum_{p=1}^M |\tilde{s}_p|-\sum_{q=1}^N |\tilde{s}'_q|\Big)$$
the numbers  $M$ (respectively, $N$)   are  the same for each of  the three chains. If not,  the following tricks help: (1) one adds and subtracts some s-patch, and thus increase $N$ and $M$ by one, and (2) for a given number $T\in \mathbb{Z}$ one takes each of the sections $T$ times and then divide  the expression by $T$.

Moreover, for the sake of consistency we should make these numbers one and the same for all the edges of the complex $K_B$, which is possible by the same reason.

Arrange the summands that appear in decompositions of $\Sigma_{ij}$, $\Sigma_{jk}$, and $-\Sigma_{ik}$ , positive ones with positive ones, and negative ones with negative ones,  to create closed triples of paths.
Apply Lemma \ref{LemmaHarmSpatch}  to each of the triples.\qed

\bigskip

Treating in similar way the rest of the $2$-faces,
we arrive at an average  of partial sections over  $Skel_2(K_B)$.
\medskip

We proceed stepwise in the similar manner. Here is how the general step looks like:

\medskip

\textbf{
Step k.}

Extend each of the partial sections  to  $Skel_k(K_B)$.  For each $(k+1)$-tuple of the vertices, say, $v_0,...,v_{k}$ of $\sigma^{n+1}$ consider the  closed $k$-chain  $$\sum_{i=0}^{k}(-1)^i\Sigma_{0...\hat{i}...k} \in C_{k-1}(\Gamma_{0...k}).$$
According to the previous section, we set $$\Sigma_{0...k}=\mathcal{H}\Big(\sum_{i=0}^{k}(-1)^i\Sigma_{0...\hat{i}...k}\Big)\in C_k(\Gamma_{0...k}).$$ The harmonic extension is taken with respect to the complex $\Gamma_{0...k}$.

The chain  $\Sigma_{0...k}$  is an average of chains supported by partial sections  over the simplex $[v_0...v_k]$. This is proven by arguments similar to those of Lemma \ref{lemma41}.

\bigskip

\textbf{Remark.} One and the same $k$-simplex $v_0,...,v_{k}$ serves as a face of several $(n+1)$-simplices.
 This construction  relies on the complex  $\Gamma_{0...k}$  only, hence it is independent on the choice of $\sigma^{n+1}$.

\bigskip
\begin{itemize}
  \item For the proof of Theorem \ref{ThmMain1} we proceed this way for all $k=1,2,...,n$, including $k=n$.  On the step $n$ we arrive at a closed $n$-chain which is an average
  of chains supported by partial sections  over the boundary of $\sigma ^{n+1}$.   Proposition \ref{prop-Euler} and the averaging idea  \textbf{complete the proof of Theorem \ref{ThmMain1}}.
  \item In the proof of Theorem \ref{ThmMain} the final step $n$ is different:
\end{itemize}

\medskip
\textbf{Step $n$.}

Our goal is to extend  the average partial sections  to  $Skel_n(K_B)$.
Now we  plan to  use the dimension-$n$  extension. On this step we further average some extended sections, but luckily
 some summands get cancelled.

For each $i$  set  $$\Sigma_{\hat{i}}= \sum_{j\in \{0,1,...,\hat{i},...,n+1\}} \ (-1)^{j+\frac{sign(j-i)+1}{2}}\ \Sigma_{012...\hat{i}...\hat{j}...n+1}.$$

By construction, $\Sigma_{\hat{i}}$ is a closed chain related to $\Gamma_{012...\hat{i}...n+1}$.

\begin{lema}\label{LemmaCancel}
               $$\sum_{i=0}^{n+1} \Sigma_{\hat{i}}=0.$$ Therefore,
               $$\sum_{i=0}^{n+1} \mathcal{W}(x,y,\Sigma_{\hat{i}})=0 \ \ \ \hbox{ for all} \ \ x,y \notin Skel_{n-1}(\Gamma_{01...n+1}).\qed$$

\end{lema}

Fix $\Sigma_{\hat{i}}$.
  For each $V_j$, $j\neq i$, we take  the associated dimension-$n$ extensions
$\mathcal{W}(V_j,\Sigma_{\hat{i}})$. Altogether we have  $(n+1)$ extensions of $\Sigma_{\hat{i}}$, each of them is an average of some s-patches of the chain $\Sigma_{\hat{i}}$ by Lemma \ref{LemmaAver}.

Therefore the chain

$$ \frac{1}{(n+1)} \sum_{j \neq i} \mathcal{W}(V_j,\Sigma_{\hat{i}})$$

is also an average of some s-patches of the chain $\Sigma_{\hat{i}}$.

The further summation by $i$  yields a closed chain which is the average of some partial sections defined on $\partial \sigma^{n+1}$:

$$\mathfrak{E}= \frac{1}{(n+1)}\sum_i \sum_{j \neq i} \mathcal{W}(V_j,\Sigma_{\hat{i}}).$$

Let us evaluate this chain at a point $x$.

$$\mathfrak{E}(x)= \frac{1}{(n+1)}\sum_i \sum_{j \neq i} \mathcal{W}(V_j,x,\Sigma_{\hat{i}})= $$

$$ \frac{1}{(n+1)}\sum_i \sum_{j} \mathcal{W}(V_j,x,\Sigma_{\hat{i}})- \frac{1}{(n+1)}\sum_i  \mathcal{W}(V_i,x,\Sigma_{\hat{i}}).  $$

Substitute $x=V_k$ and observe that the first sum vanishes by Lemma \ref{LemmaCancel}.
$$\mathfrak{E}(V_k)=- \frac{1}{(n+1)}\sum_i  \mathcal{W}(V_i,V_k,\Sigma_{\hat{i}})=\frac{1}{(n+1)}\sum_{i } \mathcal{W}(V_k,V_i,\Sigma_{\hat{i}}).$$

So $\mathfrak{E}$ is a closed $n$-chain, which is,  by construction and by Lemmata \ref{LemmaDecompN-ext} and \ref{LemmaAver}  the average of  chains supported by some partial sections defined on $\partial \sigma^{n+1}$ whose restriction on $v_i$ are $V_i$.
One more averaging over all possible choices of $V_0,...,V_{n+1}$, $V_i\in Vert(\Gamma_i)$ yields a cochain $\mathcal{E}(\sigma^{n+1})$ representing the Euler class. The proof of Theorem \ref{ThmMain} is completed. \qed

\section{Some concluding remarks}\label{sectionconcl}
\subsection*{Case $n=1$. Circle bundles.}

\begin{figure}[h]
\centering \includegraphics[width=8 cm]{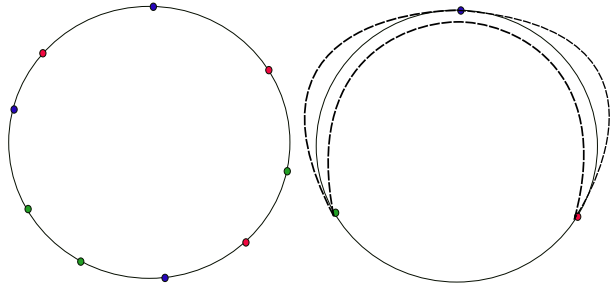}
\caption{A necklace representing $\Gamma_{012}$.  One of the eight partial sections related to a multicolored triple. This particular section contributes zero to the cochain count since the associated degree  of $S^1\rightarrow S^1$ is zero. }\label{fig3}
\end{figure}

Let us show that for $n=1$ Theorem \ref{ThmMain}  leads to the formula from \cite{M1}.
For a simplex $ \sigma^{n+1}=\sigma^2$ in the base, the complex $\Gamma_{012}$  is a triangulated circle which is viewed as a necklace with beads colored by three colors, say, red, blue, and green. The algorithm has two steps. Step $0$
fixes a multicolored triple  of beads. Combinatorially, there are two types of red-blue-green triples, positively and negatively oriented ones. Step $1$ is also the Step $n$, so in this case we do not use harmonic extensions. We only use the dimension-$n$ extensions. The chain $\Sigma_{\hat{0}}$ is a $0$-chain represented by a pair of beads $V_1$ and $V_2$.  It can be patched by one of the arches of the circle connecting $V_1$ and $V_2$ (surely, there exist other patches that wind many times around the circle, but according to our formula we use only these ones). Its dimension-$n$ extension  equals the average of these two patches. One concludes that for each triple we average over eight sections; Figure \ref{fig3} represents one of them. A simple analysis shows that each positively oriented triple of beads contributes $\frac{-1}{2}$, whereas a negative orientation  yields $\frac{1}{2}$. Thus the cochain representing the Euler class
assigns to each simplex $\sigma^2$  from the base the number $$\frac{\sharp(neg)-\sharp(pos)}{2\cdot \sharp(red)\cdot \sharp(blue)\cdot \sharp(green)}.$$

Here  $\sharp(pos)$ (respectively, $\sharp(neg)$) is the number of positively (resp., negatively) oriented multicolored triples of beads; $\sharp(red)$  denotes the number of red beads, etc.

\subsection*{Other local combinatorial formulae}

This approach leaves much room for deducing other local combinatorial formulae.

 For instance, one can use harmonic extensions for all the steps. This leads to the closest relative of N. Mn\"{e}v's formula from \cite{M2}.

  Another example. Instead of harmonic extension, one can use \textit{minimal extension:} on the first step one connects $V_i$ and $V_j$ by the shortest path in the $1$-skeleton of the complex $\Gamma_{ij}$.
By length of a path in a graph we mean the number of edges.  If there are several shortest paths, one needs to average.

On the second step one takes the patch of the "minimal area", that is, the subcomplex of the $2$-skeleton of  $\Gamma_{ijk}$  with the minimal number of $2$-cells, etc.

Another way to obtain a local combinatorial formula is to use a \textit{deformed Laplacian} instead of the usual one.
That is, one can fix an other scalar product on $C_k(\Gamma)$, for instance, to set the length of a vector associated to a cell
to be dependant on the combinatorics of the cell (say, on the number of vertices of the cell).

One expects a formula for some other cochain (representing the same cohomology class).

\end{document}